\documentclass[fleqn]{mat01}
\usepackage{times,mathtimy,amssymb,amsbsy,latexsym}
\begin{document}

\setcounter{page}{399}
\firstpage{399}

\newtheorem{theo}{Theorem}
\renewcommand\thetheo{\arabic{section}.\arabic{theo}}
\newtheorem{theor}[theo]{\bf Theorem}
\newtheorem{lem}[theo]{Lemma}
\newtheorem{pot}[theo]{Proof of Theorem}
\newtheorem{propo}[theo]{\rm PROPOSITION}
\newtheorem{rema}[theo]{Remark}
\newtheorem{defini}[theo]{\rm DEFINITION}
\newtheorem{examp}{Example}
\newtheorem{coro}[theo]{\rm COROLLARY}
\newtheorem{quest}[theo]{Question}
\def\conjecture{\trivlist\item[\hskip\labelsep{\it Conjecture.}]}
\def\exammp{\trivlist\item[\hskip\labelsep{\it Example.}]}

\font\zz=msam10 at 10pt
\def\blacksquare{\mbox{\zz{\char'245}}}

\font\xxxx=tibi at 10.4pt
\def\C{\mbox{\xxxx{C}}}

\title{Derivations into duals of ideals of Banach algebras}

\markboth{M~E~Gorgi and T~Yazdanpanah}{Ideals of Banach algebras}

\author{M~E~GORGI and T~YAZDANPANAH$^{*}$}

\address{Department of Mathematics, Faculty of Sciences, Semnan
University, Semnan, Iran\\
\noindent $^{*}$Department of Mathematics, Teacher Training University,
599, Taleghani Avenue, Tehran 15614, Iran\\
\noindent E-mail: madjideg@walla.com; 7804902@saba.tmu.ac.ir or
yazdanpanah@pgu.ac.ir}

\volume{114}

\mon{November}

\parts{4}

\Date{MS received 1 March 2004; revised 11 October 2004}

\begin{abstract}
We introduce two notions of amenability for a Banach algebra $\cal A$.
Let $I$ be a closed two-sided ideal in $\cal A$, we say $\cal A$ is
$I$-weakly amenable if the first cohomology group of $\cal A$ with
coefficients in the dual space $I^*$ is zero; i.e., $H^1({\cal
A},I^*)=\{0\}$, and, $\cal A$ is ideally amenable if $\cal A$ is
$I$-weakly amenable for every closed two-sided ideal $I$ in $\cal A$. We
relate these concepts to weak amenability of Banach algebras. We also
show that ideal amenability is different from amenability and weak
amenability. We study the $I$-weak amenability of a Banach algebra $\cal
A$ for some special closed two-sided\break ideal $I$.
\end{abstract}


\maketitle

\section{Introduction}

Let $\cal A$ be a Banach algebra and $X$ be a Banach $\cal A$-bimodule,
that is $X$ is a Banach space and an $\cal A$-bimodule such that the
module operations $(a,x)\longmapsto ax$ and $(a,x)\longmapsto xa$ from
${\cal A}\times X$ into $X$ are jointly continuous. Then $X^*$ is also a
Banach $\cal A$-bimodule if we\break define
\begin{equation*}
\langle x,ax^*\rangle =\langle xa,x^* \rangle;\quad \langle
x,x^*a\rangle = \langle ax,x^*\rangle\quad (a\in {\cal A}, ~x\in
X,~x^*\in X^*).
\end{equation*}

$\left.\right.$\vspace{-1.5pc}

We say that $Y$ is a dual $\cal A$-bimodule if there is a Banach $\cal
A$-bimodule $X$ such that $Y$ is isometrically module isomorphic with
$X^*$. Thus in particular $I$ is a Banach $\cal A$-bimodule and $I^*$ is
a dual $\cal A$-bimodule for every closed two-sided ideal $I$ in $\cal
A$.

If $X$ is  a Banach $\cal A$-bimodule, then a derivation from
$\cal A$ into $X$ is a continuous linear operator $D$ with
\begin{equation*}
D(ab)=a\cdot D(b)+D(a)\cdot b\quad (a,b\in {\cal A}).
\end{equation*}
If $x\in X$ and we define $\delta_x$ by
\begin{equation*}
\delta_x(a)=a\cdot x-x\cdot a,\quad (a\in {\cal A})
\end{equation*}
then $\delta_x$ is a derivation, and such derivations are called inner.
The space of characters on an algebra $\cal A$ is denoted by $\Phi_{\cal
A}$. Let $\varphi\in \Phi_{\cal A}\cup \{0\}$. Then $\Bbb C$ is a
symmetric $\cal A$-bimodule for the products $~a\cdot z=z\cdot
a=\varphi(a)z, (a\in {\cal A} ,z\in \Bbb C)$. In this case the bimodule
is denoted by $\Bbb C_\varphi$. A derivation from $\cal A$ into $\Bbb
C_\varphi$ is a linear functional $d$ on $\cal A$ such that
\begin{equation*}
d(ab)=\varphi(a)\cdot d(b)+d(a)\cdot \varphi(b) \quad (a,b\in {\cal A}).
\end{equation*}
Such a linear functional is called a point derivation at $\varphi$. Let
$\cal A$ be a Banach algebra. Then $\cal A$ is amenable if, whenever $D$
is a derivation from $\cal A$ to a dual $\cal A$-bimodule, then $D$ is
inner; this definition was introduced by Johnson \cite{Jo1}. $\cal
A$ is weakly amenable if, whenever $D$ is a derivation from $\cal A$ to
${\cal A}^*$, then $D$ is inner. Bade, Curtis and Dales \cite{B-C-D} have
introduced the concept of weak amenability for commutative Banach\break
algebras.

For example, it was shown in \cite{Jo1} that the group algebra $L^1(G)$
is amenable if and only if $G$ is an amenable group and in \cite{Jo2}
(or \cite{D-Gh}) that $L^1(G)$ is weakly amenable for each locally
compact group $G$.

The following definition describes the main new property that we shall
study.

\begin{defini}$\left.\right.$\vspace{.5pc}

\noindent {\rm Let $\cal A$\ be a Banach algebra and $I$ be a closed
two-sided ideal in $\cal A$. Then $\cal A$ is $I$-weakly amenable if
$H^1({\cal A},I^*)=\{0\}$; $\cal A$ is ideally amenable if $\cal A$ is
$I$-weakly amenable for every closed two-sided ideal $I$ in $\cal A$.

We begin with the following trivial  observations:
\begin{enumerate}
\renewcommand\labelenumi{(\roman{enumi})}
\leftskip .1pc
\item An amenable Banach algebra is ideally  amenable.

\item An ideally amenable Banach algebra is  weakly amenable.
\end{enumerate}}
\end{defini}

Let ${\cal A}^{\#}$ be the unitization of the commutative Banach algebra
${\cal A}$. Then for each closed two-sided ideal $I$ of ${\cal A}$
consider the following short exact sequence.
\begin{equation*}
(\Sigma)\hbox{:}\ 0\longrightarrow K\stackrel{\imath}{\longrightarrow}  {\cal
A}^{\#}\hat\otimes I \stackrel{\pi}{\longrightarrow} I\longrightarrow 0,
\end{equation*}
where $\pi$ is given by $\pi(a\otimes i)=ai$ for all $a\in{\cal
A}^{\#}$, $i\in I$, $\imath$ is the embedding map and $K= {\rm ker}\ 
\pi$. It is well-known that $B({\cal A}^{\#},I^*)$ is isometrically
isomorphic to $({\cal A}^{\#}\hat\otimes I )^*$, so we get the following
short exact sequence of linear maps:
\begin{equation*}
(\Sigma^*)\hbox{:}\ 0\longrightarrow I^*\stackrel{\pi^*}{\longrightarrow}
B({\cal A}^{\#},I^*) \stackrel{\imath^*}{\longrightarrow}
K^*\longrightarrow 0.
\end{equation*}
Let $[K;{\cal A}]={\rm Span}\,\{u\cdot a-a\cdot u\hbox{:}\ u\in K,~a\in
{\cal A}\}$, then by (\cite{Gr1}, Proposition~3.1) there is no non-zero
bounded derivation into $I^*$ if and only if $[K;{\cal A}]^{^-}=K$.

Let $\cal A$ be a Banach algebra, and let $I$ be a closed two-sided
ideal in $\cal A$. We consider the following complex
\begin{equation*}
\cdots \longrightarrow {\cal A}\hat\otimes{\cal A}\hat\otimes
I\stackrel{d_2 }{\longrightarrow} {\cal A}\hat\otimes I \stackrel
{d_1}{\longrightarrow} I \longrightarrow 0,
\end{equation*}
where the maps $d_1$ and $d_2$ are specified by the formulae:
\begin{align*}
&d_1(a\otimes i)=ai-ia\quad (a\in {\cal A}, \ i\in I);\\[.2pc]
&d_2(a\otimes b\otimes i)=b\otimes ia -ab\otimes i +a\otimes bi \
\ \ (a,b\in {\cal A} , \ i\in I).
\end{align*}
By (\cite{Jo1}, Corollary~1.3) and (\cite{Kh}, II.5.29), $H^1({\cal
A},I^*)=\{0\}$ if and only if ${\rm im}\ d_1$ is closed in $I$ and ${\rm
im}\ d_2$ is dense in ${\rm ker}\ d_1$.

\begin{theor}[\!]
Let $\cal A$ be a commutative Banach algebra and $I$ be a closed
two-sided ideal in $\cal A$. Then the following assertions are
equivalent{\rm :}
\begin{enumerate}
\renewcommand\labelenumi{\rm (\roman{enumi})}
\leftskip .35pc
\item $\cal A$ is $I$-weakly amenable.

\item ${\rm Span}\ \{(a\otimes b-b\otimes a)\cdot i-ab\otimes i{\rm :}\
a,b\in {\cal A},\ i\in I\}$ is dense in ${\cal A}\otimes I$.

\item $[K;{\cal A}]^{^-}=K$.\vspace{-1pc}
\end{enumerate}
\end{theor}

As in (\cite{B-C-D}, Theorem~1.5), we have the following theorem.

\begin{theor}[\!]
A commutative Banach algebra is weakly amenable if and only if every
derivation from $\cal A$ into a commutative Banach $\cal A$-bimodule is
zero.
\end{theor}

By Theorem~1.3 we conclude that a commutative Banach algebra $\cal A$ is
weakly amenable if and only if it is ideally amenable. As in
(\cite{B-C-D}, Theorem~3.14), let $K$ be an infinite, compact metric
space and let $\alpha\in (0,1/2)$ then ${\rm lip}_\alpha K$ is weakly
amenable but it is not amenable, therefore this is an example of ideally
amenable Banach algebra, that is not amenable.

We now consider an example of non-commutative Banach algebra that is
ideally amenable but it is not amenable.

\begin{exammp}
Let ${\cal A}=\ell^1(\Bbb N)$. We define the product on ${\cal A}$ by
$f\cdot g=f(1)g$ for all $f$ and $g$ in ${\cal A}$. It is obvious that
${\cal A}$ is a Banach algebra with this product and norm
$\parallel\cdot\parallel _1$. It is straightforward that ${\cal A}$ has
no approximate identity, so by (\cite{Jo1}, Proposition~1.6) ${\cal A}$
is not amenable. Let $I$ be a closed two-sided ideal of ${\cal A}$, it
is easy to see that if $I\neq {\cal A}$, then $I\subseteq \{f\in {\cal
A};~f(1)=0\}$.

Let $I\neq{\cal A}$ and $D\hbox{:}\ {\cal A}\longrightarrow I^*$ be a
derivation. For $f\in {\cal A }$ we consider
\begin{align*}
&\tilde{f}\hbox{:}\ \Bbb N\longrightarrow \Bbb C,\\[.2pc]
&\tilde{f}(n) = \left\{ \begin{array}{ccc} 
f(n), &&n\geq 2\\[.2pc]
0, &&n=1
\end{array} \right.,
\end{align*}
then $f=f\cdot e+\tilde f$ such that
\begin{equation*}
e(n) = \left\{\begin{array}{ccc}
1, &&n=1\\[.2pc]
0, && n\neq1
\end{array} \right..
\end{equation*}
Therefore $D(f)=f(1)D(e)+D(\tilde f)$ but we have $\tilde D(\tilde
f\cdot e)=\tilde f\cdot D(e)+ D(\tilde f)\cdot e$, $\tilde f \cdot
e=\tilde f\cdot D(e)=0$ and $D(\tilde f)\cdot e=D(\tilde f)$.
Consequently $D(\tilde f)=0$. Let $x^*=-D(e)$.
\begin{align*}
f\cdot x^*-x^*\cdot f &= 0-f(1)x^*\\[.2pc]
&= -f(1)(-D(e))\\[.2pc]
&= f(1)D(e)\\[.2pc]
&= D(f).
\end{align*}
So $D=\delta_{x^*}$ and $H'({\cal A},I^*)=\{0\}$. If $I=\cal A$,
then by Assertion 1 of \cite{Zh}, $H^{1}({\cal
A},I^*)=\{0\}$.\break \phantom{0} \hfill $\blacksquare$\vspace{.7pc}
\end{exammp}

\looseness 1 Let $\cal A$ be a Banach algebra with a bounded right(left) approximate
identity and let $X$ be a Banach $\cal A$-bimodule on which $\cal A$
acts trivially on the left(right). Then by (\cite{Jo1}, Proposition~1.5)
$H^{1}({\cal A},X^*)=\{0\}$, so if $I$ is a closed two-sided ideal in
${\cal A}$ and ${\cal A} I= \{0\} (I {\cal A}=\{0\})$, then ${\cal A}$
is $I$-weakly amenable. If $G$ is a \hbox{commutative}\break \newpage \noindent non-discrete group, then
by [Gh-L-W], $L^1(G)^{**}$ is not weakly amenable. Let $J=\{F\in
L^1(G)^{**}\hbox{:}\ L^1(G)^{**}F=0\}$. Then $J$ is a closed two-sided
ideal in $L^1(G)^{**}$ and $L^1(G)^{**}$ is $J$-weakly amenable.

\begin{theor}[\!] 
Let $\cal A$ be a Banach algebra and $I$ be a closed two-sided ideal in
$\cal A$ and $\cal A$ be $I$-weakly amenable{\rm ,} let $\varphi\in
\Phi_{\cal A}${\rm ,} such that $I\not\subseteq\ {\rm ker}~\varphi$. Then
there is no non-zero point derivation at $\varphi$.
\end{theor}

\begin{proof}
Let $d\hbox{:}\ {\cal A}\longrightarrow \Bbb C_\varphi$ be a non-zero
point derivation, and let $\pi\hbox{:}\ {\cal A}^*\longrightarrow I^*$
be the adjoint of $\imath\hbox{:}\ {I}\longrightarrow {\cal A}$.
Consider the map $D\hbox{:}\ {\cal A} \longrightarrow I^*$ defined by
$D(a)=d(a)\pi(\varphi)$. It is easy to see that $D$ is a derivation.
Since $\cal A$ is $I$-weakly amenable, there exists $\lambda\in I^* $
with $D(a)=a\cdot\lambda-\lambda \cdot a~~(a\in {\cal A})$. Take $i\in
I$ with $\varphi(i)=1$. If $\hbox{ ker} \ \varphi \subseteq\ \hbox{ker}\
d$, then there exists $\alpha \in \Bbb C$ such that $d= \alpha \varphi$.
So $2\alpha = 2 \alpha \varphi(i)=2 d(i)= 2 d(i) \varphi(i) = d(i^2) =
\alpha \varphi(i^2) =\alpha$. Therefore, $\alpha=0$ and this is a
contradiction. Consequently $\hbox{ker}\ \varphi\!\nsubseteq\!\hbox{ker}\
d$ and there exists $a\in \ \hbox{ker}\ \varphi$ with $d(a)=1$. Set
$i'=i+(1-d(i))ia=i+ia -d(i)ia$. Then $\varphi(i')=d(i')=1$ and so
\begin{equation*}
1=(D(i'))(i')=(i'\cdot\lambda-\lambda\cdot i')(i')= \lambda(i'{^2}) -
\lambda(i'{^2})=0,
\end{equation*}
is a contradiction.\hfill $\blacksquare$
\end{proof}

\begin{theor}[\!]
Let ${\cal A}$ be a weakly amenable Banach algebra and for each closed
two-sided ideal I such that $I=\overline{{\cal A}I\bigcup I{\cal A}}${\rm ,}
${\cal A}$ is $I$-weakly amenable. Then ${\cal A}$ is ideally amenable.
\end{theor}

\begin{proof}
Let $I$ be a closed two-sided ideal in ${\cal A}$. Put
$J=\overline{{\cal A}I \bigcup I{\cal A}}$. It is easy to see that $J$
is a closed two-sided ideal in ${\cal A}$ and $J=\overline{J{\cal A}
\bigcup {\cal A} J}$. Let $\imath\hbox{:}\ J \longrightarrow I$ be the
natural embedding and $D\hbox{:}\ {\cal A} \longrightarrow I^*$ be a
derivation. Then $\imath^*\circ D\hbox{:}\ {\cal A} \longrightarrow J^*$
is a derivation. So there is a $m \in J^*$ such that $\imath^*\circ D
=\delta_m$. Let $x^* $ be the extension of $m$ into $I$ by Hahn--Banach
theorem. For every $a, b \in {\cal A}$ we have
\begin{align*}
\langle i , D(ab) \rangle &= \langle ia , D(b) \rangle + \langle bi , D(a) \rangle\\[.2pc]
&= \langle \imath(ia), D(b) \rangle + \langle \imath(bi ) , D(a) \rangle \\[.2pc]
&= \langle ia, bm - mb \rangle + \langle bi, am - ma \rangle\\[.2pc]
&= \langle iab - bia +bia -abi, m \rangle \\[.2pc]
&= \langle i, abx^* -x^*ab \rangle = \langle i, \delta_{x^*} (ab) \rangle ~~~~(i \in I).
\end{align*}
Hence $D(ab) = \delta_{x^*} (ab)$. Since ${\cal A}$ is weakly
amenable, so ${\cal A}^2 = {\cal A}$ and therefore
$D=\delta_{x^*}$ and $D$ is inner.\hfill $\blacksquare$
\end{proof}

\begin{theor}[\!]
Let $\cal A$ be a Banach algebra{\rm ,} $X$ be a Banach $\cal
A$-bimodule and $Y$ be a closed submodule of $X$. If $H^1({\cal
A},Y^*)=\{0\}$ and $H^1({\cal A},({X/Y})^*)=\{0\}${\rm ,} then
$H^1({\cal A},X^*)=\{0\}$.
\end{theor}

\begin{proof}
Let $D\hbox{:}\ {\cal A}\longrightarrow X^*$ be a derivation, and
$\pi\hbox{:}\ X^*\longrightarrow Y^*$ be adjoint of $\imath\hbox{:}\
Y\longrightarrow X$. Then $\pi$ is a $\cal A$-bimodule homomorphism and
$\pi\circ D\hbox{:}\ {\cal A} \longrightarrow Y^*$ is a derivation.
Therefore there exists $y^*\in Y^*$ such that $\pi\circ D=\delta_{y^*}$.
Let $x^*$ be an extension of $y^*$ by Hahn--Banach theorem. Then
$d=D-\delta_{x^*}\hbox{:}\ {\cal A}\longrightarrow Y^\perp$ is a
derivation, but $Y^\perp=(X/Y)^*$. Therefore there exists $x^*_1\in
Y^\perp\subseteq X^*$ such that $d=\delta_{x^*_1}$, so
$D=\delta_{(x^*+x^*_1)}$.\hfill $\blacksquare$
\end{proof}

\begin{coro}$\left.\right.$\vspace{.5pc}

\noindent Let $I$ and $J$ be closed two-sided ideals in Banach algebra
$\cal A${\rm ,} and $J\subseteq I$. If $\cal A$ is $J$\!-weakly amenable
and $H^1({\cal A},(I/J)^*)=\{0\}${\rm ,} then $\cal A$ is $I$-weakly
amenable.
\end{coro}

\begin{coro}$\left.\right.$\vspace{.5pc}

\noindent Let $I$ be a closed two-sided ideal in Banach algebra $\cal
A$. If $H^1({\cal A},({\cal A}/I)^*)=\{0\}$ and $\cal A$ is $I$-weakly
amenable{\rm ,} then $\cal A$ is weakly amenable.
\end{coro}

\begin{theor}[\!]
Let ${\cal A}$ be a Banach algebra and $I$ be a closed two-sided ideal
in ${\cal A}$ with a bounded approximate identity. If ${\cal A}$ is
ideally amenable{\rm ,} then $I$ is ideally amenable.
\end{theor}

\begin{proof}
Let $J$ be a closed two-sided ideal in $I$. It is easy to see that $J$
is an ideal in ${\cal A}$. Let $D\hbox{:}\ I \longrightarrow J^*$ be a
derivation. By (\cite{Ru}, Proposition~2.1.6), $D$ can be extended to a
derivation ${ \tilde D}\hbox{:}\ {\cal A} \longrightarrow J^*$. So there
is a $m \in J^*$ such that ${\tilde D}= \delta_{m}$. Then $D(i)={\tilde
D}(i)=\delta_m$ for each $i \in I$. So $D$ is inner.\hfill
$\blacksquare$
\end{proof}

\begin{coro}$\left.\right.$\vspace{.5pc}

\noindent Let ${\cal A}$ be a Banach algebra with a bounded approximate
identity and ${\cal M({\cal A})}$ be the multiplier algebra of ${\cal A}$. If
${\cal M({\cal A})}$ is ideally amenable{\rm ,} then ${\cal A}$ is ideally
amenable.
\end{coro}

\begin{defini}{\rm \cite{Gr3}}$\left.\right.$\vspace{.5pc}

\noindent {\rm Let ${\cal A}$ be a Banach algebra and $I$ be a closed
two-sided ideal in ${\cal A}$. We say that $I$ has the trace extension
property, if every $ m \in I^*$ such that $ am = ma $ for each $a \in
{\cal A} $, can be extended to $ a^* \in A^*$ such that $ aa^* = a^* a$
for each $ a \in {\cal A}$.}
\end{defini}

Let $I$ be a closed two-sided ideal in Banach algebra $\cal A$.
Johnson in \cite{Jo1} has shown that, if ${\cal A}$ is amenable,
then $I$ is amenable if and only if $I$ has a bounded approximate
identity, and $\cal A$ is amenable if $I$ and ${\cal A}/I$ are amenable.
Also if ${\cal A}/I$ and $I$ are weakly amenable, then $\cal A$ is
weakly amenable \cite{Gr2}. Also if $I$ has the trace extension property
and $\cal A$ is weakly amenable, then ${\cal A}/I$ is weakly amenable
\cite{Gr2}. We prove a similar proposition for ideal amenability.

\begin{theor}[\!]
Let ${\cal A}$ be an ideally amenable Banach algebra and $I$ be a closed
two-sided ideal in ${\cal A}$ that has the trace extension property.
Then{\rm ,} ${\cal A} /I$ is ideally amenable.
\end{theor}

\begin{proof}
Let $J/I$ be a closed two-sided ideal in ${\cal A}/ I$. Then $J$ is a
closed two-sided ideal in ${\cal A}$. We write $\pi\hbox{:}\ J
\longrightarrow J/I$, $q\hbox{:}\ {\cal A} \longrightarrow {\cal A} /I$
for the natural quotient maps and $\pi ^*$ for the adjoint of $\pi$. Let
$D\hbox{:}\ {\cal A}/I \longrightarrow (J/I)^*$ be a derivation. Then
$d= \pi^* \circ D \circ q\hbox{:}\ {\cal A} \longrightarrow J^*$ is a
derivation, so there exists $x^* \in J^*$ such that $d = \delta _{x^*}$.
Let $m$ be the restriction of $x^*$ to $I$. Then $m \in I^*$ and for all
$i \in I$ we have
\begin{align*}
\langle i, am -ma \rangle &= \langle ia -ai , m \rangle = \langle ia
-ai, x^*\rangle \\[.2pc]
&= \langle i , \delta _{x^*} (a) \rangle = \langle i , \pi^* \circ D
\circ q (a) \rangle \\[.2pc]
&= \langle \pi (i) , D \circ q (a) \rangle = \langle I , D(a+I)\rangle\\[.2pc]
&= 0.
\end{align*}
Therefore, $am = ma (a \in {\cal A})$, and so $a^* \in {\cal A} ^*$ such
that $ aa^* = a^* a$ and $a^*$ is an extension of $m$. Let $y^*$ be
the restriction of $ a^* $ to $J$. Then $ y^* \in J^*$ and $ x^* - y^*
=0 $ on $I$. Therefore, $x^* - y^* \in (J/I) ^* $ and we have
\begin{align*}
\langle j+ I, D( a+ I)\rangle &= \langle \pi (j), D(q(a))\rangle \\[.2pc]
&= \langle j, \pi^* \circ D \circ q (a)\rangle \\[.2pc]
&= \langle j , \delta _ {x^*} (a) \rangle = \langle j , \delta _{x^* -
y^* } (a) \rangle\\[.2pc]
&= \langle j + I , \delta _{x^* - y^* } (a+I).
\end{align*}
Hence $D= \delta _{x^* - y^*}$ and therefore ${\cal A}/I$ is ideally
amenable.\hfill $\blacksquare$
\end{proof}

Let $I$ be a closed two-sided ideal in Banach algebra $\cal A$. We
consider some easy remarks about the relations between $I$-weak
amenability of $\cal A$ and weak amenability of $I$ and $\cal A$. 

As in \cite{Jo2} we know that the group algebra $L^1(G)$ is weakly
amenable for every locally compact group $G$, but we do not know whether
or not $L^1(G)$ is ideally amenable. By the following theorem we can
show that $M(G)$ is $I$-weakly amenable if and only if $L^1(G)$ is
$I$-weakly amenable for every closed two-sided ideal $I$ in $L^1(G)$.

\begin{theor}[\!]
Let $\cal A$ be a Banach algebra and let $J$ be a closed two-sided ideal
in $\cal A$ with a bounded approximate identity. Then for every closed
two-sided ideal $I$ in $J${\rm ,} $J$ is $I$-weakly amenable if and only if
$\cal A$ is $I$-weakly amenable.
\end{theor}

\begin{proof}
Let $(j_\alpha)_{\alpha\in \Lambda}$ be a bounded approximate identity
for $J$ and $D\hbox{:}\ J \longrightarrow I^*$ be a derivation. Consider
the map $\tilde D\hbox{:}\ {\cal A}\longrightarrow I^*$ defined by
\begin{equation*}
\tilde D(a)=w^*-\lim_\alpha (D(aj_\alpha)-a\cdot D(j_\alpha))\quad (a\in
\cal A).
\end{equation*}
By (\cite{Ru}, Proposition~2.1.16) $\tilde D$ is a continuous
derivation. If $D=0$, then $\tilde D=0$, since $JI=IJ=I$. Therefore
$H^1(J,I^*)=H^1({\cal A},I^*)$, and this implies that $\cal A$ is
$I$-weakly amenable if and only if $J$ is $I$-weakly amenable.\hfill
$\blacksquare$
\end{proof}

Let ${\cal A}^{\#} $ be the unitization of ${\cal A}$. We know that
${\cal A}$ is amenable if and only if ${\cal A}^{\#}$ is amenable. If
${\cal A}$ is weakly amenable then ${\cal A}^{\#}$ is weakly amenable
(see (\cite{D-Gh-G}, Proposition~1.4, (ii))). For ideal amenability we
have the following:

\begin{propo}$\left.\right.$\vspace{.5pc}

\noindent Let ${\cal A}$ be a Banach algebra. Then ${\cal A}$ is ideally
amenable if and only if ${\cal A}^{\#}$ is ideally amenable.
\end{propo}

\begin{proof}
Let ${\cal A}^{\#}$ be ideally amenable, $I$ be a closed two-sided ideal
of ${\cal A}$ and $D\hbox{:}\ {\cal A} \longrightarrow I^*$ be a
derivation. It is easy to see that $I$ is a closed two-sided ideal of
${\cal A} ^{\#}$, and ${\tilde D}\hbox{:}\ {\cal A}^{\#} \longrightarrow
I^*$ such that $ {\tilde D} (a+\alpha) = D(a)$, $ (a \in {\cal A} ,~
\alpha \in \Bbb C )$ is a derivation. Since ${\cal A}^{\#}$ is ideally
amenable, $ \tilde{D} $ is inner and hence $D$ is inner.

Conversely, let ${\cal A}\ $ be ideally amenable and $I$ be a closed
two-sided ideal in ${\cal A}^{\#}$. Since ${\cal A}$ is ideally
amenable, therefore ${\cal A}$ is weakly amenable, so ${\cal A}^{\#}$ is
weakly amenable. Therefore, we can suppose that $1 \notin I$ and $I$ is
an ideal of ${\cal A}$. Let $D\hbox{:}\ {\cal A}^{\#} \longrightarrow
I^*$ be a derivation, then $ D(1)=0$ and we can consider $D$ as a
derivation from ${\cal A}$ into $I^*$ and therefore $D$ is\break
inner.\hfill $\blacksquare$
\end{proof}

Let $\cal A$ be a non-unital Banach algebra. Then similar to the
Proposition~1.14 we can show that ${\cal A}^{\#}$ is not $\cal A$-weakly
amenable. Let $\cal A$ be the augmentation ideal of $L^1(PLS_2(R))$.
Michael White has shown that $\cal A$ is not weakly amenable and ${\cal
A}^{\#}$ is an example of weakly amenable Banach algebra that is not
ideally amenable.

\section{${\pmb{\cal C}}^*$-algebras}

Recall first that, a ${\cal C}^*$-algebra is amenable if and only if it
is nuclear \cite{Ha}. However, every ${\cal C}^*$-algebra is weakly
amenable \cite{Ha}. We show that every ${\cal C}^*$-algebra is ideally\break
amenable.

\setcounter{theo}{0}
\begin{lem}
Let ${\cal A}$ be a Banach algebra and $I$ be a closed two-sided ideal
in ${\cal A}$. If $I$ is weakly amenable{\rm ,} then ${\cal A}$ is
$I$-weakly amenable.
\end{lem}

\begin{proof}
Let $D\hbox{:}\ {\cal A} \longrightarrow I^*$ be a derivation and
$\imath\hbox{:}\ I \longrightarrow {\cal A} $ be the embedding map. Then
$D\circ \imath\hbox{:}\ I \longrightarrow I^*$ is a derivation, and so
there exists $ m \in I^*$ such that $ D \circ \imath = \delta_m $. Since
$I$ is weakly amenable, $\overline{I^2}=I$, and for $i,j \in I$ we have
\begin{align*}
\langle ij, D(a)\rangle &= \langle i, jD(a)\rangle = \langle i,
D(ja)-D(j)a \rangle \\[.2pc]
&= \langle i, jam - mja \rangle -\langle ai, jm -mj \rangle \\[.2pc]
&= \langle ija, m \rangle - \langle aij, m \rangle = \langle ij, am -
ma \rangle \\[.2pc]
&= \langle ij, \delta_{m} (a) \rangle\quad (a \in {\cal A)}.
\end{align*}

Therefore $D=\delta_m$, and so $D$ is inner.\hfill $\blacksquare$
\end{proof}

\begin{coro}$\left.\right.$\vspace{.5pc}

\noindent Every $C^*$-algebra is ideally amenable.
\end{coro}

\begin{proof}
Let ${\cal A}$\ be a $C^*$-algebra and $I$ be a two-sided closed ideal in
${\cal A}$, then $I$ is a $C^*$-algebra by (\cite{B-D}, Theorem~38.18),
so $I$ is weakly amenable. By Lemma~2.1, ${\cal A}$ is $I$-weakly
amenable.\hfill $\blacksquare$
\end{proof}

\begin{coro}$\left.\right.$\vspace{.5pc}

\noindent Let $H$ be a Hilbert space{\rm ,} $K(H)$ be the space of
compact operators in $B(H)$. Then $H^1(B(H),K(H)^*)=\{0\}$.
\end{coro}

\begin{coro}$\left.\right.$\vspace{.5pc}

\noindent Let $H$ be an infinite dimensional Hilbert space. Then $B(H)$
is ideally amenable{\rm ,} but not amenable.
\end{coro}

Let $G$ be a locally compact topological group. Then $L^1(G)$ is an
ideal in $L^1(G)^{**}$ if and only if $G$ is compact, then by Lemma~2.1
we have the following result.

\begin{coro}$\left.\right.$\vspace{.5pc}

\noindent Let $G$ be a compact topological group. Then $L^1(G)^{**}$ is
$L^1(G)$-weakly amenable.
\end{coro}

\section{Maximal ideals}

Let $M$ be a maximal ideal in $\cal A$. We study conditions when
$\cal A$ is $M$-weakly amenable.

\setcounter{theo}{0}
\begin{theor}[\!]
Let $\cal A$ be a unital weakly amenable Banach algebra. Then for every
maximal ideal $ M$ in $\cal A${\rm ,} $\cal A$ is $ M$-weakly amenable.
\end{theor}

\begin{proof}
Let $D\hbox{:}\ {\cal A}\longrightarrow { M}^*$ be a derivation. For
each $a$ in $\cal A$, $D(a)$ has an extension $\tilde{D}(a)\in \cal
A^*$ such that $\tilde{D}(1)=0$. For $ a,b \in \cal A$
there exists $\lambda_1 , \lambda_2 \in \Bbb{C}$ and $m,n \in M$ such
that
\begin{align*}
\tilde{{D}}(ab) &= \tilde {D}((\lambda_1\cdot 1+m)(\lambda_2\cdot 1+n))\\[.2pc]
&= \lambda_1 {D}(n)+\lambda_2 {D}(m) +\tilde {D}(mn)\\[.2pc]
&= \lambda_1 {D}(n)+\lambda_2 {D}(m) +m\cdot {D}(n)+
{D}(m)\cdot n\\[.2pc]
&= (\lambda_1\cdot 1+m)\cdot {D}(n)+{D}(m)\cdot(\lambda_2+n)\\[.2pc]
&= a\cdot\tilde {D}(b)+\tilde {D}(a)\cdot b.
\end{align*}
Therefore ${\tilde D}\in {\cal Z}^1 ({\cal A},{\cal A}^*)$ and there
exists $a^*\in {\cal A}^*$ such that ${\tilde D}={\delta}_{a^*}$.
Obviously ${D}=\delta_{a^*\mid M}$.\hfill $\blacksquare$
\end{proof}

\begin{theor}[\!]
Let $\cal A$ be a Banach algebra and $M$ be a closed two-sided ideal in
$\cal A$ of codimension one. If $\cal A$ is $M$-weakly amenable{\rm ,}
then $M$ is weakly amenable.
\end{theor}

\begin{proof}
Let ${D}\hbox{:}\ M\longrightarrow M^*$ be a derivation. Since ${\cal
A}=M\oplus \Bbb C$, the map ${D}_1\hbox{:}\ {\cal A}\longrightarrow M^*$
defined by ${ D}_1(m+c)={ D}(m)$ for $m\in {\cal M}, \alpha \in \Bbb C$
is an inner derivation. Consequently ${D}$ is inner.\hfill $\blacksquare$
\end{proof}

Let $F_2$ be the free group on two generators. Let $\ell^0(F_2)=\{\mu\in
{\ell^1(F_2)}\hbox{:}\ \mu(F_2)=0\}$. Then by Theorem~3.1, $\ell^1(F_2)$
is $\ell^0(F_2)$-weakly amenable, and by Theorem~3.2, $\ell^0(F_2)$ is
weakly amenable.

\begin{lem}
Let $\cal A$ be a unital commutative Banach algebra and $M$ be a closed
two-sided ideal in $\cal A$ of codimension one. If $\cal A$ is
$M$\!-weakly amenable{\rm ,} then $\cal A $ is weakly amenable.
\end{lem}

\begin{proof}
By the above theorem, we know that $M$ is weakly amenable. On the other
hand, we have ${\cal A}=M\oplus\Bbb C$. Therefore by (\cite{Gr1},
Proposition~2.3), $\cal A$ is weakly amenable.\hfill $\blacksquare$
\end{proof}

Now we have the following theorem:

\begin{theor}[\!]
Let $\cal A$ be a commutative unital Banach algebra{\rm ,} the following
assertions are equivalent{\rm :}
\begin{enumerate}
\renewcommand\labelenumi{\rm (\roman{enumi})}
\leftskip .35pc
\item $\cal A$ is weakly amenable{\rm ,}

\item $\cal A$ is ideal weakly amenable{\rm ,}

\item $\cal A$ is $M$-weakly amenable for some maximal ideal $M$
in $\cal A${\rm ,}

\item $H^1({\cal A},X^*)=\{0\}$ for every commutative
Banach $\cal A$-bimodule $X$.
\end{enumerate}
\end{theor}

\begin{propo}$\left.\right.$\vspace{.5pc}

\noindent Let $\cal A$ be a commutative Banach algebra with a bounded
approximate identity and $M$ be a maximal modular ideal in $\cal A$ also
${\cal A}_z\cap M=0$ {\rm (}where ${\cal A}_z$ is the set of topological
divisor of zero elements in $\cal A${\rm )}. Let $D{\rm :}\ {\cal
A}\longrightarrow M$ be a derivation such that $D(M)=\{0\}${\rm ,} then\break 
$D=0$.
\end{propo}

\begin{proof}
Let $(e_\alpha)_{\alpha\in I}$ be a bounded approximate identity for
$\cal A$. We may suppose that $e_\alpha\in {\cal A}\setminus M$ for
every $\alpha\in I$. Let $a\in\cal A$ and $D(a)\neq 0$. Then for every
$\alpha\in I$ there exists $a_\alpha\in {\cal A}\setminus M$ and
$m_\alpha\in M$ such that
\begin{equation}
e_\alpha=m_\alpha+a_\alpha a.
\end{equation}
Therefore $0=D^2(e_\alpha)=2D(a_\alpha)D(a)$. Since $D(a)\in M$ and
${\cal A}_z\cap M=0$, for every $\alpha \in I$, $D(a)=0$. On the other
hand, for every $a'\in\cal A$, we have
\begin{align*}
D(a') &= \lim_\alpha D(e_\alpha a')\\[.2pc]
&= \lim_\alpha D(e_\alpha)a' +\lim_\alpha e_\alpha D(a')\\[.2pc]
&= \lim_\alpha D(e_\alpha)a'+D(a').
\end{align*}
Consequently $\displaystyle\lim_\alpha D(e_\alpha)a'=0$ for every $a'\in
A$. Let $a'$ be a non-zero element of $M$, since ${\cal A}_z\cap M=0$,
$\displaystyle\lim_\alpha D(e_\alpha)=0$ and by (1), we have
$\displaystyle\lim_\alpha a_\alpha D(a)=0$. Therefore
$\displaystyle\lim_\alpha a_\alpha=0$. Now for each $b \in \cal A$,
$b=\displaystyle\lim_\alpha be_\alpha = \displaystyle\lim_\alpha
b(m_\alpha + a_\alpha a)= \displaystyle\lim_\alpha [b(m_\alpha) +b
a_\alpha a]= \displaystyle\lim_\alpha bm_\alpha$. Consequently
$(m_\alpha)_{\alpha\in I}$ is a approximate identity for $\cal A$ but by
(\cite{Pa}, Theorem~5.2.7) $M$ is closed and so $M=\cal A$ and this is a
contradiction.\hfill $\blacksquare$
\end{proof}

\section{Problems}

We are interested in the problems listed below.

Johnson \cite{Jo1} has shown that $\cal A \widehat{\otimes} {\cal B}$
is amenable whenever $\cal A$ and $\cal B$ are amenable Banach algebras.
So we can raise the following question.

\setcounter{theo}{0}
\begin{quest}
{\rm If $\cal A$ and $\cal B$ are ideally amenable Banach algebras, then is
$\cal A \widehat{\otimes}\cal B$ ideally amenable?}
\end{quest}

We know that $L^1(G)$ is amenable if and only if $G$ is an amenable
group \cite{Jo1}, and also $L^1(G)$ is weakly amenable for every locally
compact group (\cite{Jo1} or \cite{D-Gh}).

\begin{quest}{\rm 
Under what conditions the group algebra $L^1(G)$ is ideally amenable?}
\end{quest}

\begin{quest}{\rm 
If $L^1(G)^{**}$ is ideally amenable, then so are $M(G)$ and $L^1(G)$
ideally amenable?}
\end{quest}

\begin{quest}{\rm 
Is $\ell^1(F_2)$ ideally amenable?}
\end{quest}

For a Banach algebra $\cal A$, the amenability of $\cal A^{**}$
necessitates the amenability of $\cal A$ (\cite{Da}, Proposition 2.8.59) and similarly for weak
amenability provided $\cal A$ is a left ideal in $\cal A^{**}$
\cite{Gh-L-W}. So we can raise the following question.

\begin{quest}{\rm 
If $\cal A^{**}$ is ideally amenable, then is $\cal A$ ideally amenable?}
\end{quest}
\pagebreak

\section*{Acknowledgement}

We are grateful to Prof.~Alireza Medghalchi for his kind
encouragement and invaluable suggestions. Also we would like to thank
the referee for some interesting and useful suggestions which improved
the paper.


\begin{thebibliography}{9999999}
\bibitem[B-C-D]{B-C-D} Bade~W~G, Curtis~P~G and Dales~H~G, Amenability and weak
amenability for Beurling and Lipschitz algebra, {\it Proc. London Math.
Soc. (3)} {\bf 55} (1987) 359--377 

\bibitem[B-D]{B-D} Bonsall~F~F and Duncan~J, Complete normed algebras
(Berlin: Springer) (1973) 

\bibitem[Da]{Da} Dales H~E, Banach algebra and automatic continuity, London Mathematical Society Monographs (Oxford: Clarendon Press) (2000) vol.~24

\bibitem[D-Gh-G]{D-Gh-G} Dales~H~G, Ghahramani~F and Gronbaek~N, Derivations
into iterated duals of Banach algebras, {\it Studia Math.} {\bf 128}
(1998) 19--54 

\bibitem[D-Gh]{D-Gh} Despic~M and Ghahramani~F, Weak amenability of group
algebras of locally compact groups, {\it Canad. Math. Bull.} {\bf 37}
(1994) 165--167 

\bibitem[Gh-L-W]{Gh-L-W} Ghahramani~F, Loy~R~J and Willis~G~A, Amenability and
weak amenability of second conjugate Banach algebras, {\it Proc. Am.
Math. Soc.} {\bf 124(5)} (1996) 1489--1497 

\bibitem[Gr1]{Gr1} Gronbak~N, A characterization of weakly amenable Banach
algebras, {\it Studia Math.} {\bf 94} (1989) 150--162 

\bibitem[Gr2]{Gr2} Gronbak~N, Weak amenability of group algebras, {\it Bull.
London Math. Soc.} {\bf 23} (1991) 231--284 

\bibitem[Gr3]{Gr3} Gronbak~N, Weak and cyclic amenable for non-commutative
Banach algebras, {\it Proc. Edinburg Math. Soc.} {\bf 35} (1992)
315--328 

\bibitem[Ha]{Ha} Haagerup~U, All nuclear ${\cal C}^*$-algebras are amenable,
{\it Invent. Math.} {\bf 74} (1983) 305--319 

\bibitem[Jo1]{Jo1} Johnson~B~E, Cohomology in Banach algebras, {\it Mem. Am.
Math. Soc.} {\bf 127} (1972) 

\bibitem[Jo2]{Jo2} Johnson~B~E, Weak amenability of group algebras, {\it
Bull. London Math. Soc.} {\bf 23} (1991) 281--284 

\bibitem[Kh]{Kh} Khelmeskii~A~Ya, The homology of Banach and topological
algebras (Kluwer) (1989) 

\bibitem[Pa]{Pa} Palmer~Theodore~W, Banach algebra and the general theory of
$*$-algebras (Cambridge University Press) (1994) vol.~1 

\bibitem[Ru]{Ru} Runde~V, Lectures on amenability (Berlin, Heidelberg, New
York: Springer-Verlag) (2001) 

\bibitem[Zh]{Zh} Zhang~Yong, Weak amenability of a class of Banach algebras,
{\it Canad. Math. Bull.} {\bf 44} (2001) 504--508
\end{thebibliography}
\end{document}